\documentclass[11pt]{article}
\usepackage{amssymb}
\usepackage{amsmath}
\usepackage{amsthm}
\usepackage{latexsym}
\usepackage{amsfonts}
\usepackage{graphicx}
\usepackage{graphics}

\theoremstyle{plain}
\newtheorem{thm}{Theorem}[section]   

\newtheorem{lem}[thm]{Lemma}

\theoremstyle{definition}

\newtheorem*{Proof}{Proof}

\newcommand{\bi} {{\beta}}

\newcommand{\ld} {{\ldots}}
\newcommand{\sm} {{\smallsetminus}}

\newcommand{\de} {{\delta}}
\newcommand{\De} {{\varDelta}}
\newcommand{\si} {{\sigma}}
\newcommand{\la} {{\lambda}}
\newcommand{\el} {{\ell}}
\newcommand{\e} {{\varepsilon}}

\newcommand{\mi} {{\mu}}

\newcommand{\dis}{\displaystyle}
\newcommand{\ssum}{\sum\limits}

\newcommand{\cp}{{\cal{P}}}
\newcommand{\ch}{{\cal{H}}}

\newcommand{\cn}{{\cal{N}}}

\newcommand{\ra}{{\rightarrow}}

\newcommand{\qb}{$\quad\blacksquare$}
\def\1{\it1\hspace*{-0.150cm}{\footnotesize{I}}}

\def\R{{\mathbb{R}}}
\def\C{{\mathbb{C}}}
\def\Q{{\mathbb{Q}}}
\def\N{{\mathbb{N}}}

\begin{document}
\title{\bf Common hypercyclic vectors for families of backward shift operators}
\author{\bf N. Tsirivas}
\footnotetext{{The research project is implemented within the framework of the Action ``Supporting Postdoctoral Researchers'' of the Operational Program ``Educational and Lifelong Learning'' (Action's Beneficiary: General Secretariat for Research and Technology), and is co-financed by the European Social Fund (ESF) and the Greek State.}}

\date{}
\maketitle
{\bf Abstract:} We provide necessary and sufficient conditions on the existence of common hypercyclic vectors for multiples of the backward shift operator along sparse powers. Our main result strongly generalizes corresponding results which concern the full orbit of the backward shift. Some of our results are valid in a more general context, in the sense that they apply for a wide class of hypercyclic operators.    

%
\noindent


\section{Introduction}\label{sec1}
\noindent

We consider the space $\el^2$ of square summable sequences over the field of complex numbers $\mathbb{C}$ endowed with the topology that is induced by the $\el^2$ norm $\|\cdot\|_2:\el^2\ra\R^+$, where
\[
\|x\|_2:=\bigg(\sum^{+\infty}_{j=1}|x_j|^2\bigg)^{1/2} \ \ \text{for every} \ \ x=(x_1,x_2,\ld,x_n,\ld)\in\el^2.
\]
We write $\|\cdot\|:=\|\cdot\|_2$ for simplicity. Let $B$ be the unweighted backward shift operator on $\el^2$, that is
\[
B(x_1,x_2,x_3,\ld)=(x_2,x_3,\ld), \ \ \text{for} \ \ (x_1,x_2,\ld)\in\el^2.
\]
Let $\la\in\C$, $|\la |>1$ and consider the set of hypercyclic vectors for $\la B$, that is
\[
\ch C(\la B):=\big\{x=(x_1,x_2,\ld)\in\el^2\big|\overline{\big\{(\la B)^n
(x),\;n=1,2,\ld\big\}}=\el^2\big\}.
\]
A comprehensive treatment on hypercyclicity can be found in the books \cite{2}, \cite{10}. For the reader's convenience we include the relevant definition. A sequence of continuous operators $(T_n)$ acting on a Frechet space $X$ is called hypercyclic provided there exists a vector $x\in X$ so that the set $\{ T_n(x): n=1,2,\ldots \}$ is dense in $X$. Such a vector is called hypercyclic for $(T_n)$ and the set of hypercyclic vectors for $(T_n)$ is denoted by $\ch C(\{ T_n\} )$. When $T_n$ comes from the iterates of a single operator we sat that $T$ is hypercyclic and $\ch C(T)$ denotes the set of hypercyclic vectors for $T$, i.e. $$\ch C (T)=\{ x\in X: \{ T^nx: n=1,2,\ldots \} \,\,\textrm{is dense in} \,\, X\}.$$    
It is well known that for every $\la \in \mathbb{C}$ with $|\la |>1$ the set $\ch C(\la B)$ is a dense, $G_\de$ subset of $(\el^2,\|\cdot\|)$ and as the reader may guess Baire's theorem should be involved in the arguments. The following question arises naturally.
If we fix an uncountable subset $J\subset\{z\in\C|\,|z|>1\}$ is it true that
$\bigcap\limits_{\la\in J}\ch C(\la B)\neq\emptyset$\,\,?
In this direction, Abakumov and Gordon \cite{AbGo} proved that:
\[
\bigcap_{|\la |>1}\ch C(\la B)\neq\emptyset, 
\]
the best possible result one can expect concerning the existence of common hypercyclic vectors for multiples of the backward shift. Later on, Costakis and Sambarino \cite{CoSa} gave a different proof of this result, which, roughly speaking, is based on the so called common hypercyclicity criterion. In this criterion, Baire's category theorem appears. Actually, Costakis and Sambarino showed that $\bigcap\limits_{|\la |>1}\ch C(\la B)$ is a $G_\de$ and dense subset of $(\el^2,\|\cdot\|)$; hence non-empty. What is interesting here is the uncountable range of $\la$'s, which makes things harder if one wishes to apply Baire's theorem.   

One can refine the above problem as follows.
Let $(k_n)$ be a fixed subsequence of natural numbers. It is known, and very easy to prove, that the sequence $( (\la B)^{k_n})$ is also hypercyclic, that is, there exists $x\in \el^2$ such that the set $\{ (\la B)^{k_n}(x):n=1,2,\ldots \}$ is dense in $\el^2$. Such a vector is called hypercyclic for $( (\la B)^{k_n})$ and the set of these vectors is denoted by $\ch C(\{(\la B)^{k_n}\})$. From the above it should be also clear, or at least expected, that $\ch C(\{(\la B)^{k_n}\})$ is $G_\de$ and dense subset of $(\el^2,\|\cdot\|)$.  

Now we are ready to ask the following 

{\bf Question:}
\textit{Fix a strictly increasing sequence $(k_n)$ of natural numbers. For which uncountable sets $J\subset \{ \lambda \in \mathbb{C} :|\la |>1 \}$ }, 

\[
\bigcap\limits_{\la\in J}\ch C(\{(\la B)^{k_n}\})\neq\emptyset\,\,\text{?}
\]

It turns out that the answer to this question depends heavily on the sequence $(k_n)$. In particular, what matters is how sparse the sequence $(k_n)$ has been chosen. Our main result is the following
\begin{thm} \label{mainresult}
Let $(k_n)$ be a strictly increasing sequence of positive integers. 
\begin{enumerate}
\item[(i)] If $\ssum^{+\infty}_{n=1}\dfrac{1}{k_n}<+\infty$ then $\bigcap\limits_{\la\in I}\ch C(\{(\la B)^{k_n}\})=\emptyset$ for every non-degenerate interval in $(1,+\infty )$.
\item[(ii)] If $\ssum^{+\infty}_{n=1}\dfrac{1}{k_n}=+\infty$ then the set $\bigcap\limits_{\la\in (1,+\infty)}\ch C(\{(\la B)^{k_n}\})$ is residual in $\el_2$, i.e. it contains a $G_{\de}$ and dense set in $\el_2$; hence $\bigcap\limits_{\la\in (1,+\infty)}\ch C(\{(\la B)^{k_n}\})\neq \emptyset$.
\item[(iii)] If $\ssum^{+\infty}_{n=1}\dfrac{1}{k_n}=+\infty$ there exists a $G_{\delta }$ and dense subset $\cp$ in $\{ \la \in \mathbb{C}: |\la |>1 \}$ with full $2$-dimensional Lebesgue measure in $ \{ \lambda \in \mathbb{C}: |\la |>1 \}$ such that $\bigcap_{\la\in\cp}\ch C(\{(\la B)^{k_n}\}) $ is  residual in $\el_2$. In particular, $\bigcap_{\la\in\cp}\ch C(\{(\la B)^{k_n}\})\neq\emptyset$.
\end{enumerate}
\end{thm}

Unfortunately we are unable to show whether $\cp $ in item $(iii)$ of Theorem \ref{mainresult} can be replaced by $\{ \la\in \mathbb{C}: |\la |>1 \}$. So, this remains an open problem. On the other hand, both items $(i)$ and $(iii)$ hold in a more general setting (there is nothing special if one choses to work with the backward shift) and this is evident if one follows the relevant proofs, see sections $2$ and $4$. For instance, the interested readers will have no difficulties in formulating general statements for items $(i)$ and $(iii)$ that involve operators $T$ so that for a given sequence of positive integers $(k_n)$, the sequence$((\lambda T)^{k_n})$ is hypercyclic for every $\lambda$ lying in some interval or annulus, possibly with infinite length or infinite area. We mention that a kind of similar line of research is pursued in \cite{Ba}, \cite{8}, \cite{Tsi}, \cite{Tsi2}, \cite{Tsi3}, where questions similar to the above one are studied for translation type operators acting on the space of entire functions. Results on the existence of common hypercyclic vectors for uncountable families of operators and, in particular, of backward shift operators can be found in \cite{AbGo}, \cite{3}, \cite{ChaSa2}, \cite{ChaSa4}, \cite{5}, \cite{CoSa}, \cite{LeMu}, \cite{San}, \cite{13}.  
Our paper is organized as follows. Each one of the following sections $2$, $3$, $4$, is devoted to the proof of items $(i)$, $(ii)$, $(iii)$ of Theorem \ref{mainresult} respectively. 

The proof of item $(i)$ relies on an estimate which concerns the size (in terms of Lebesgue measure) of the set$\{z\in\C|\;|z^{n}\cdot w-1|<\e \} \cap [m, M]$, for given $w\in \mathbb{C}$, $\e >0$, $1<m<M$. This approach is implicit in \cite{2}, \cite{3}, \cite{13} and refines an idea of Borichev. The common hypercyclicity criterion due to Costakis and Sambarino
cannot be applied in order to conclude item $(ii)$. What we do, is to refine in a sense this criterion in the particular case of backward shift. It seems plausible that our method will possibly work for other operators as well. We mention that there are quite a few, relatively new and powerful, criteria establishing the existence of common hypercyclic vectors for uncountable families of operators, see \cite{2}, \cite{3}, \cite{13}. However, it is not clear to us whether these criteria can be used in our case. Finally, the proof of item $(iii)$relies on the following three ingredients: $1)$ item $(ii)$, $2)$ a metric result of Weyl which says that, if $(k_n)$ is a given sequence of distinct integers then the sequence $(k_nx)$ is uniformly distributed mod $1$, see Theorem 4.1 in \cite{KN}, and $3)$ Cavalieri's principle, see page 149 in \cite{Hal}. Actually, to prove item $(iii)$ we elaborate on the proof of Proposition 5.2 from \cite{BaCo}.

\section{A negative result}\label{sec2}
\noindent

Fix a subsequence $(k_n)$ of natural numbers such that $\ssum^{+\infty}_{n=1}\dfrac{1}{k_n}<+\infty$. Our goal, in this section, is to show that $\bigcap\limits_{\la\in J}\ch C(\{(\la B)^{k_n}\})=\emptyset$.
First of all we need two helpful lemmas.
\begin{lem}\label{lem2.1}
Let $z_0\in\C$, $N_0\in\N$, $\e_0\in(0,1)$ be three fixed numbers.

We consider the set
\[
G:=G(z_0,N_0,\e_0)=\{z\in\C|\;|z^{N_0}\cdot z_0-1|<\e_0\}.
\]
Then $G$ is an open subset of $\C$.
\end{lem}
The proof of Lemma \ref{lem2.1} is straightforward and is omitted. We proceed with the following

\begin{lem}\label{lem2.2}
Let $z_0\in\C$, $N_0\in\N$, $N_0>1$, $\e_0\in(0,1)$, $\mi_0>1$, $M_0>\mi_0$ be fixed. Using the notation of the previous lemma, the following estimate holds.
\[
\la(G(z_0,N_0,\e_0)\cap[\mi_0,M_0])\le M_0\cdot\bigg(\sqrt[N_0]{\dfrac{1+\e_0}{1-\e_0}}-1\bigg).
\]
\end{lem}
\begin{Proof}
Suppose that $G_1:=G(z_0,N_0,\e_0)\cap[\mi_0,M_0]\neq\emptyset$, (the other case is trivial for the proof). Since $G_1$ is open it contains two different elements $a_0$ and $A_0$ where $a_0<A_0$. By the definition of $G(z_0,N_0, \e_0)$ it follows that
\setcounter{equation}{0}
\begin{eqnarray}
|a^{N_0}_0z_0-1|<\e_0.  \label{eq1}
\end{eqnarray}
We have $z_0\neq0$. Using (\ref{eq1}) and the triangle inequality we get
\begin{eqnarray}
|z_0|>\frac{1-\e_0}{a_0^{N_0}}.  \label{eq2}
\end{eqnarray}
Let $z_0:=t_1+t_2i$, where $t_1=Re(z_0)$ and $t_2=Im(z_0)$. Then,
\begin{eqnarray}
|a^{N_0}_0z_0-1|^2=|z_0|^2a^{2N_0}_0+1-2t_1a^{N_0}_0.  \label{eq3}
\end{eqnarray}
By (\ref{eq1}) and (\ref{eq3}) we conclude that
\begin{eqnarray}
|z_0|^2(a^{N_0}_0)^2-2t_1(a^{N_0}_0)+(1-\e^2_0)<0.  \label{eq4}
\end{eqnarray}
Consider the trinomial
\[
A(y):=|z_0|^2y^2-2t_1y+(1-\e^2_0).
\]
Because $A(a^{N_0}_0)<0$, by (\ref{eq4}) and $|z_0|^2>0$ $(z_0\neq0)$ the trinomial $A(y)$ has positive discriminant $\De>0$ and two roots $\rho_1$ and $\rho_2$, where $\rho_1<\rho_2$ and $A^{N_0}_0<\rho_2$ $(\ast)$ (it follows easily by the above that $\rho_2>0$ as well). Hence,
\begin{eqnarray}
\rho_2\le\frac{1+\e_0}{|z_0|}.  \label{eq5}
\end{eqnarray}
By (\ref{eq2}), $(\ast)$ and (\ref{eq5}) we get
\begin{eqnarray}
A_0^{N_0}<\frac{1+\e_0}{1-\e_0}\cdot a^{N_0}_0\Rightarrow A_0-a_0<M_0\cdot
\bigg(\sqrt[N_0]{\frac{1+\e_0}{1-\e_0}}-1\bigg).  \label{eq6}
\end{eqnarray}
So, for every $a,A\in G_1$, $a<A$,
\[
A-a<M_0\cdot\bigg(\sqrt[N_0]{\frac{1+\e_0}{1-\e_0}}-1\bigg).  \eqno{(6)'}
\]
This gives that
\begin{eqnarray}
\de(G_1)\le M_0\bigg(\sqrt[N_0]{\frac{1+\e_0}{1-\e_0}}-1\bigg), \ \ \text{where} \ \ \de(G_1) \ \ \text{is the diameter of} \ \ G_1  \label{eq7}
\end{eqnarray}
The set $G_1$ is open in $[\mi_0,M_0]$ and bounded, so
\begin{eqnarray}
G_1\subseteq[\inf G_1,\sup G_1]\Rightarrow\la(G_1)\le\la([\inf G_1,\sup G_1])=\sup G_1-\inf G_1=\de(G_1).  \label{eq8}
\end{eqnarray}
By (\ref{eq7}) and (\ref{eq8}) we arrive at 
\[
\la(G_1)\le M_0\cdot\bigg(\sqrt[N_0]{\frac{1+\e_0}{1-\e_0}}-1\bigg)
\]
and the proof of Lemma \ref{lem2.2} is complete. \qb
\end{Proof}

Now we are ready to prove the main theorem of this section.
\begin{thm}\label{thm2.1}
Let $(k_n)$ be a subsequence of natural numbers such that $\ssum^{+\infty}_{n=1}\dfrac{1}{k_n}<+\infty$. Then
\[
\bigcap_{\la\in(1,+\infty)}\ch C(\{(\la B)^{k_n}\})=\emptyset.
\]
\end{thm}
\begin{Proof}
We fix two positive numbers $\mi_0,M_0$ such that $1<\mi_0<M_0<+\infty$. It suffices to prove that
\[
\bigcap_{\la\in[\mi_0,M_0]}\ch C(\{(\la B)^{k_n}\})=\emptyset .
\]

We set
\[
\si_0:=\sum^{+\infty}_{n=1}\frac{1}{k_n}<+\infty 
\]
and fix some positive number $\de_0\in\Big(0,\dfrac{M_0-\mi_0}{\si_0}\Big)$. For instance we may take $\de_0:=\dfrac{M_0-\mi_0}{2\si_0}$. It is obvious that $\dfrac{\de_0}{M_0}>0$, so $e^{\frac{\de_0}{M_0}}>1$. Because $\dis\lim_{x\ra0^+}\dfrac{1+x}{1-x}=1$, there exists some positive number $\e_0\in(0,1)$ such that
\setcounter{equation}{0}
\begin{eqnarray}
\frac{1+\e_0}{1-\e_0}<e^{\frac{\de_0}{M_0}}.  \label{eq1}
\end{eqnarray}
We now fix some positive number $\e_0\in(0,1)$ such that (\ref{eq1}) holds.

%
%
There exists some natural number $N_0\in\N$ such that
\begin{eqnarray}
\bigg(\frac{\de_0}{M_0N}+1\bigg)^N>\frac{1+\e_0}{1-\e_0}  \label{eq2}
\end{eqnarray}
for every $N\in\N$, $N\ge N_0$, by (\ref{eq1}).
To arrive at a contradiction, suppose that $\bigcap\limits_{\la\in[\mi_0,M_0]}\ch C(\{(\la B)^{k_n}\})\neq\emptyset$. We fix some
\[
x_0:=(x_1,x_2,\ld)\in\bigcap_{\la\in[\mi_0,M_0]}\ch C(\{(\la B)^{k_n}\})
\]
and let $e_1:=(1,0,0,\ld)\in\el^2$. We fix some $\la_0\in[\mi_0,M_0]$. Since
\[
x_0\in\bigcap_{\la\in[\mi_0,M_0]}\ch C(\{(\la B)^{k_n}\}),
\]
there exists a subsequence $(\mi_n)$ of $(k_n)$ such that $(\la_0B)^{\mi_n}(x_0)\ra e_1$.
For $\e_1:=\e_0/2$, there exists some natural number $n_0\ge N_0$ such that
\begin{eqnarray}
\|(\la_0B)^{\mi_n}(x_0)-e_1\|<\e_1 \ \ \text{for every} \ \ n\ge n_0.  \label{eq3}
\end{eqnarray}
By (\ref{eq3}) we get
\begin{eqnarray}
|\la_0^{\mi_n}x_{\mi_n+1}-1|<\e_1 \ \ \text{for every} \ \ n\in\N, \ \ n\ge n_0.  \label{eq4}
\end{eqnarray}
Though $\la_0$ is fixed, we apply the above for every $\la\in [\mi_0,M_0]$. Hence, (\ref{eq4}) implies that for every $\la\in[\mi_0,M_0]$ there exists some natural number $v\ge N_0$ such that
\begin{eqnarray}
|\la^{k_v}x_{k_v+1}-1|<\e_1.  \label{eq5}
\end{eqnarray}

Setting
\[
L:=\big\{\la\in[\mi_0,M_0]/\,\exists\;v\ge N_0:|\la^{k_v}x_{k_v+1}-1|<\e_1\big\},
\]
we get $L=[\mi_0,M_0]$, by (\ref{eq5}).
Consider the set
\[
\cn_1:=\big\{v\in\N/v\ge N_0 \;\text{and}\;\exists\;\la\in[\mi_0,M_0]\;\text{such that}\;|\la^{k_v}x_{k_v+1}-1|<\e_1\big\}.
\]
Then, $\cn_1\neq\emptyset$ by (\ref{eq5}).

Let $v\in\cn_1$ and define the set
\[
G_v:=\big\{\la\in[\mi_0,M_0]\bigg|\;|\la^{k_v}x_{k_v+1}-1|<\e_1\big\}.
\]
It is obvious that $G_v\neq\emptyset$ for every $v\in\cn_1$.
By Lemma \ref{lem2.1}, $G_v$ is open in $[\mi_0,M_0]$ for every $v\in\cn_1$ and 
\begin{eqnarray}
[\mi_0,M_0]=\bigcup_{v\in\cn_1}G_v.  \label{eq6}
\end{eqnarray}
By the properties of Lebesgue measure, Lemma \ref{lem2.2} and (\ref{eq6}) we have
\begin{align}
M_0-\mi_0&=\la([\mi_0,M_0])=\la\bigg(\bigcup_{v\in\cn_1}G_v\bigg)\le\sum_{v\in\cn_1}
\la(G_v)\nonumber \\
&\le\sum_{v\in\cn_1}M_0\bigg(\sqrt[k_v]{\frac{1+\e_1}{1-\e_1}}-1\bigg) \nonumber\\
&\le M_0\cdot\sum^{+\infty}_{v=N_0}\bigg(\sqrt[k_v]{\frac{1+\e_1}{1-\e_1}}-1\bigg). \label{eq7}
\end{align}
Observe that $\dfrac{1+\e_0}{1-\e_0}>\dfrac{1+\e_1}{1-\e_1}$, since
\begin{eqnarray}
\e_1\in(0,\e_0).  \label{eq8}
\end{eqnarray}
By (\ref{eq2}) and (\ref{eq8}) we get
\begin{eqnarray}
\sqrt[k_v]{\frac{1+\e_1}{1-\e_1}}-1<\frac{\de_0}{M_0k_v} \ \ \text{for every}\ \ v\ge N_0  \label{eq9}
\end{eqnarray}
and (\ref{eq7}), (\ref{eq9}) imply that
\begin{eqnarray}
M_0\cdot\sum^{+\infty}_{v=N_0}\bigg(\sqrt[k_v]{\frac{1+\e_1}{1-\e_1}}-1\bigg)<
\de_0\cdot\sum^{+\infty}_{v=N_0}\frac{1}{k_v}<M_0-\mi_0  \label{eq10}
\end{eqnarray}
(by the definition of $\de_0$). Obviously, the inequalities (\ref{eq7}) and (\ref{eq10}) are in contradiction and this completes the proof of Theorem \ref{thm2.1}. \qb
\end{Proof}
\section{The positive case in the half line $(1,+\infty )$} \label{sec3}
\noindent

Throughout this section  we fix a subsequence $(k_n)$ of natural numbers such that $\ssum^{+\infty}_{n=1}\dfrac{1}{k_n}=+\infty$. We shall prove the following
\begin{thm}\label{thm3.1}
The set $\bigcap_{\la\in(1,+\infty)}\ch C(\{(\la B)^{k_n}\})$ is a residual subset of $(\el^2,\|\cdot\|)$.
\end{thm}

In order to prove Theorem \ref{thm3.1} we assign some notations and terminology. Let
$D:=\big\{x=(x_1,x_2,\ld)\in\el^2\big|\{n\in\N\big|x_n\neq0\}$ is a finite subset of $\N$ and $x_n\in\Q+i\Q$ for every $n=1,2,\ld\big\}$, where $\Q$ is the set of rational numbers,
The set $D$ is countable and dense in $(\el^2,\|\cdot\|_2)$. We set $\overline{0}:=(0,0,\ld)\in\el^2$ and $D^\ast:=D\sm\{\overline{0}\}$.
Let $\varPsi:=\{y_1,y_2,y_3,\ld,y_n,\ld\}$ be an enumeration of $D^\ast$.  We fix a strictly decreasing sequence of positive numbers $(a_n)$ such that $a_n\ra1$ (for example $a_n:=1+\dfrac{1}{n}$, $n=1,2,\ld\;$) and we also fix a strictly increasing sequence $(\bi_n)$ of positive numbers such that $\bi_n\ra+\infty$ and $a_1<\bi_1$ (for example $\bi_n=n+2$, $n=1,2,\ld\;$). Then we set $\De_n:=[a_n,\bi_n]$, $n=1,2,\ld\;$. Of course, the sequence of compact sets $\De_n$, $n=1,2,\ld$ forms an exhausting family of $(1,+\infty)$, i.e.  $(1,+\infty)=\bigcup\limits^{+\infty}_{n=1}\De_n$.

For every $n,j,s,m\in\N$ let us define  
\[
E_{\De_n}(j,s,m):=\{x=(x_1,x_2,\ld)\in\el^2|
\]
for every $\la\in\De_n$ there exists some $v\in\N$, $v\le m$ such that
\[
\|(\la B)^{k_v}(x)-y_j\|<\frac{1}{s}\bigg\}.
\]
We finally set
\[
G:=\bigcap^{+\infty}_{n=1}\bigcap^{+\infty}_{j=1}\bigcap^{+\infty}_{s=1}
\bigcup^{+\infty}_{m=1}E_{\De_n}(j,s,m).
\]
\begin{lem}\label{lem3.2}
For every $n,j,s,m\in\N$ the set $E_{\De_n}(j,s,m)$ is open in $(\el^2,\|\cdot\|)$.
\end{lem}
\begin{Proof}
The proof is easy and we leave it to the reader.
\end{Proof}
\begin{lem}\label{lem3.3}
We have
\[
G\subseteq\bigcap_{\la\in(1,+\infty)}\ch C(\{(\la B)^{k_n}\})
\]
\end{lem}
\begin{Proof}
The proof is trivial. \qb
\end{Proof}

We proceed with two more lemmas that lie in the heart of the argument.
\begin{lem}\label{lem3.4}
Let $a_0,b_0$ be two positive real numbers such that $1<a_0<b_0<+\infty$. Let $(k_n)$ be a subsequence of natural numbers such that $\ssum^{+\infty}_{n=1}\dfrac{1}{k_n}=+\infty$.
Then for every positive number $\e>0$ there exists a finite number of terms $k_{n_0},k_{n_0+1},\linebreak\ld,k_{n_0+i_0}$ of $(k_n)$ for some natural numbers $n_0,i_0$ and positive numbers $\bi_1,\bi_2,\ld,\bi_{i_0+1}$ such that: for every $\la\in[a_0,b_0]$, there exists some $j\in\{0,1,\ld,i_0\}$ such that
\[
|\la^{k_{n_0+j}}\bi_{j+1}-1|<\e.
\]
\end{lem}
\begin{Proof}
We fix some positive number $\e_0\in(0,1)$. After we fix some natural number $n_0$ such that:
\setcounter{equation}{0}
\begin{eqnarray}
k_{n_0}>\frac{\log(1+\e_0)}{\log\Big(\dfrac{b_0}{a_0}\Big)}.  \label{eq1}
\end{eqnarray}
Of course $\ssum^{+\infty}_{j=0}\dfrac{1}{k_{n_0+j}}=+\infty$. So by (\ref{eq1}) there exists the unique natural number $i_0\in\{0,1,2,\ld\}$ such that
\begin{eqnarray}
(1+\e_0)^{\ssum^{i_0}_{j=0}\frac{1}{k_{n_0+j}}}\le\frac{b_0}{a_0} \ \ \text{and} \label{eq2}
\end{eqnarray}
\begin{eqnarray}
(1+\e_0)^{\ssum^{i_0+1}_{j=0}\frac{1}{k_{n_0+j}}}>\frac{b_0}{a_0}.  \label{eq3}
\end{eqnarray}
We set $\bi_1:=\dfrac{1}{a^{k_{n_0}}}$. Then for every $\la\in[a_0,b_0]$ such that $\la<a_0\cdot\sqrt[k_{n_0}]{1+\e_0}$, we get $|\la^{k_{n_0}}\cdot\bi_1-1|<\e_0$.

After we set $a_1:=a_0\cdot\sqrt[k_{n_0}]{1+\e_0}$ and $\bi_2:=\dfrac{1}{a_1^{k_{n_0+1}}}$.

Then, for every $\la\in[a_1,b_0]$ with $\la<a_1\cdot\sqrt[k_{n_0+1}]{1+\e_0}$ we have $|\la^{k_{n_0+1}}\bi_2-1|<\e_0$.

We continue inductively.

We suppose that we have defined the number $a_i=a_0(1+\e_0)^{\ssum^{i-1}_{j=0}\frac{1}{k_{n_0+j}}}$ for some $i\in\{1,2,\ld,i_0-1\}$.

After we define $\bi_{i+1}:=1/a^{k_{n_0+i}}_i$ and for every $\la\in[a_i,b_0]$ with $\la<a_i\cdot\sqrt[k_{n_0+i}]{1+\e_0}$ we get $|\la^{k_{n_0+i}}\cdot\bi_{i+1}-1|<\e_0$.

Because $a_0<a_1<a_2<\cdots<a_{i_0}\le b_0$ and $a_{i_0}$ is the maximum number in $(a_0,b_0]$ that we can obtain with the above procedure, after a finite number of steps we exclude the interval $[a_0,b_0]$ and we get the conclusion of this lemma with the following data:

We have
\[
a_i=a_0\cdot\prod^{i-1}_{j=0}(1+\e_0)^{\frac{1}{k_{n_0+j}}} \ \ \text{for every} \ \ i=1,2,\ld,i_0
\]
\[
\bi_{i+1}=1/a^{k_{n_0+i}}_i \ \ \text{for every} \ \ i=0,1,\ld,i_0.
\]
It follows that for every $\la\in[a_i,b_0]$ where $\la<a_i\cdot\sqrt[k_{n_0+i}]{1+\e_0}$ we have:
\[
|\la^{k_{n_0}+i}\bi_{i+1}-1|<\e_0 \ \ \text{for every} \ \ i=0,1,\ld,i_0.
\]
With the above data the proof of this lemma is completed. \qb
\end{Proof}
\begin{lem}\label{lem3.5}
Let $(k_n)$ be a subsequence of natural numbers such that $\ssum^{+\infty}_{n=1}\dfrac{1}{k_n}=+\infty$.

Then for every positive number $M>0$ there exists a subsequence $(\mi_n)$ of $(k_n)$ such that:

(i) $\mi_{n+1}-\mi_n>M$ for every $n=1,2,\ld$ and

(ii) $\ssum^{+\infty}_{n=1}\dfrac{1}{\mi_n}=+\infty$.
\end{lem}
\begin{Proof}
We fix some positive number $M_0>1$.

Let $N_0:=[M_0]+1\ge2$ (where $[x]$ is the integer part of $x\in\R$).

We consider the subsequences $(\mi^j_\rho)_{\rho=1,2,\ld}=(k_{\rho N_0+j})_{\rho=1,2,\ld}$ of $(k_n)$ for every $j=0,1,\ld,N_0-1$, that is $\mi^j_\rho:=k_{\rho N_0+j}$, $\rho=1,2,\ld$ for $j=0,1,\ld,N_0-1$.

We fix some $j_0\in\{0,1,\ld,N_0-1\}$.  We consider the subsequence $(\mi^{j_0}_\rho)_{\rho=1,2,\ld}$ of $(k_n)$. \vspace*{0.2cm} \\
\noindent
{\bf Claim 1.}

For every $v_1,v_2\in\N$, $v_1<v_2$ we have: $k_{v_2}-k_{v_1}\ge v_2-v_1$. \vspace*{0.2cm} \\
\noindent
{\bf Proof of the Claim 1.}

Because $(k_n)$ is a subsequence of natural numbers we have:
\[
\begin{array}{l}
  k_{v_1+1}-k_{v_1}\ge1, \; k_{v_1+2}-k_{v_1+1}\ge1,\ld, \\ [1.5ex]
  k_{v_1}+(v_2-v_1)-k_{v_1}+(v_2-v_1)-1 \ge1
\end{array}
\]
Adding the above inequalities we take the conclusion of Claim 1.

We apply Claim 1 for the terms of subsequence $(\mi^{j_0}_\rho)_{\rho=1,2,\ld}$ and we have:
\setcounter{equation}{0}
\begin{eqnarray}
\mi^{j_0}_{\rho+1}-\mi^{j_0}_\rho>M_0.  \label{eq1}
\end{eqnarray}
By (\ref{eq1}) we have that each one from the subsequences $(\mi^j_\rho)_{\rho=1,2,\ld}$ for $j=0,1,\ld,N_0-1$ of $(k_n)$ satisfies (\ref{eq1}). It is obvious, by the definition of the sequences $(\mi^j_\rho)_{\rho=1,2,\ld}$, for $j=0,1,\ld,N_0-1$ that these do not have common terms by pairs and every term of the sequence $(k_v)$ belongs in a unique subsequence $(\mi^j_\rho)_{\rho=1,2,\ld}$ for some $j\in\{0,1,\ld,N_0-1\}$.

This gives that
\begin{eqnarray}
\sum^{+\infty}_{n=1}\frac{1}{k_n}=\sum^{N_0-1}_{j=0}\sum^{+\infty}_{\rho=1}\frac{1}
{\mi^j_\rho}.  \label{eq2}
\end{eqnarray}
Because $\ssum^{+\infty}_{n=1}\dfrac{1}{k_n}\!=\!+\infty$, the relation (\ref{eq2}) gives us that there exists one $j\!\in\!\{0,1,\ld,N_0-1\}$ at least such that $\ssum^{+\infty}_{\rho=1}\dfrac{1}{\mi^j_\rho}=+\infty$, where the subsequence $(\mi^j_\rho)$ of $(k_n)$ satisfies the properties (i) and (ii) obviously and this completes the proof of this lemma. \qb
\end{Proof}

After the above preparation we are ready now to prove the following lemma that is the basic result that gives us Theorem \ref{thm3.1}.
\begin{lem}\label{lem3.6}
For every $n,j,s\in\N$ the set $\bigcup\limits^{+\infty}_{m=1}E_{\De_n}(j,s,m)$ is dense in $(\el^2,\|\cdot\|)$.
\end{lem}
\begin{Proof}
We fix $n_0,j_0,s_0\in\N$ and we will show that the set $\bigcup\limits^{+\infty}_{m=1}E_{\De_{n_0}}(j_0,s_0,m)$ is dense in $(\el^2,\|\cdot\|)$.

We set $E:=\bigcup\limits^{+\infty}_{m=1}E_{\De_{n_0}}(j_0,s_0,m)$ for simplicity.

Let $y_{j_0}:=(q_1,q_2,\ld,q_{v_0},0,0,\ld)$ where $q_{v_0}\neq0$, $y_{j_0}(v)=q_v$, for every $v=1,2,\ld$, $y_{j_0}(v)=0$ for every $v\ge v_0+1$, and $q_j\in\Q+i\Q$ for every $j=1,2,\ld,v_0$, for some fixed $v_0\in\N$.

We fix $\e_0>0$ and $c_0=(c_1,c_2,\ld,c_{v_1},0,0,\ld)\in D$ where $c_{v_1}\neq0$, $v_1\in\N$, fixed, $c_0(v)=c_v$ for every $v=1,2,\ld$, $c_0(v)=0$ for every $v\ge v_1+1$, $c_j\in\Q+i\Q$ for every $j=1,2,\ld,v_1$. We consider the ball $S_{\el^2}(c_0,\e_0):=\{x\in\el^2\big|\|x-c_0\|<\e_0\}$.

We will show that
\[
E\cap S_{\el^2}(c_0,\e_0)\neq\emptyset.  \eqno{(\ast)}
\]
In order to show the relation $(\ast)$ it suffices to show that there exists some $x_0=(x_1,x_2,\ld,x_n,\ld)\in\el^2$ and $m_0\in\N$ such that

(i) $\|x_0-c_0\|<\e_0$ and

(ii) for every $\la\in \De_{n_0}$ there exists some $v\in\N$, $v\le m_0$ such that
\setcounter{equation}{0}
\begin{eqnarray}
\|(\la B)^{k_v}(x_0)-y_{j_0}\|<\frac{1}{s_0}.  \label{eq1}
\end{eqnarray}
We will succeed (i) and (ii) above as follows:

From the data of the problem we define a finite number of complex numbers $x_j$, $j=1,2,\ld,\el_0$ for some fixed $\el_0\in\N$.

Afterwards, we define the sequence $x_0:=(x_1,x_2,\ld,x_{\el_0},0,0)$ where $x_0(j)=x_j$ for every $j=1,2,\ld,\el_0$ and $x_0(j)=0$ for every $j\ge\el_0+1$. So we have $x_0\in\el^2$.

We define also a natural number $m_0$.

Finally, we show that $x_0$ and $m_0$ satisfy properties (i) and (ii) of (\ref{eq1}) as above.

Without loss of generality let $\De_{n_0}:=[a_0,b_0]$, where $1<a_0<b_0<+\infty$. We set
\[
M_1:=\max\{|q_j|,\;j=1,2,\ld,v_0\}>0.
\]

We also set
\[
M_2:=\frac{1}{2\log a_0}\cdot\log\bigg(\frac{2s^2_0\cdot M^2_1\cdot v_0}{1-\dfrac{1}{a_0}}\bigg).
\]
By Lemma \ref{lem3.5} we choose a subsequence $(\mi_n)$ of $(k_n)$ such that the following two properties hold:

(i) $\mi_{n+1}-\mi_n>\max\{M_2,v_0\}$ for every $n=1,2,\ld$ and

(ii) $\ssum^{+\infty}_{n=1}\dfrac{1}{\mi_n}=+\infty$.

We set $\e_1:=\dfrac{1}{\sqrt{2v_0}\cdot s_0\cdot M_1}$.

Now we can choose some fixed natural number $v_2\in\N$ such that the following three inequalities hold:
\begin{enumerate}
\item[a)] $\mi_{v_2}>v_1+1$
\item[b)] $\mi_{v_2}>\dfrac{\log(1+\e_1)}{\log\Big(\dfrac{b_0}{a_0}\Big)}$
\item[c)] $\mi_{v_2}>\dfrac{1}{2\log a_0}\cdot\log\bigg(\dfrac{v_0\cdot M^2_1}{\e^2_0\cdot\Big(1-\dfrac{1}{a_0}\Big)}\bigg)$.
\end{enumerate}
Because $\ssum^{+\infty}_{n=1}\dfrac{1}{\mi_n}=+\infty$ we have $\ssum^{+\infty}_{j=0}\dfrac{1}{\mi_{v_2+j}}=+\infty$.

Let $i_0$ be the unique natural number $i_0\in\{0,1,2,\ld\}$ such that
\begin{eqnarray}
(1+\e_1)^{\ssum^{i_0}_{j=0}\frac{1}{\mi_{v_2+j}}}\le\dfrac{b_0}{a_0}\ \  \text{and} \ \ (1+\e_1)^{\ssum^{i_0+1}_{j=0}\frac{1}{\mi_{v_2+j}}}>\dfrac{b_0}{a_0} \label{eq2}
\end{eqnarray}
We set $m_1:=v_2+i_0$. The natural number $m_1$ is the natural number in order (\ref{eq1}) holds. That is let $m_0$ be the unique natural number such that: $k_{m_0}:=\mi_{m_1}$. Then $m_0$ is the natural number we search in (\ref{eq1}).

Now, we are ready to define the vector $x_0=(x_1,x_2,\ld)\in\el^2$ straightforward with full details.

We define $x_j=c_j$ for every $j=1,2,\ld,v_1$. We define $x_j=0$ for every $j\in\N$ such that: $v_1+1\le j\le\mi_{v_2}$.

Because $\ssum^{+\infty}_{n=1}\dfrac{1}{\mi_n}=+\infty$, we can apply Lemma \ref{lem3.4} for the sequence $(\mi_n)$. By the relations (b) and (\ref{eq2}) above that the numbers $v_2$ and $\e_1$, $i_0$ satisfy and using Lemma \ref{lem3.4} for the sequence $(\mi_n)$ we take that there exists a finite number of positive numbers $\bi_1,\bi_2,\ld,\bi_{i_0+1}$ that are defined completely in Lemma \ref{lem3.4} from our data such that for every $\la\in[a_0,b_0]$ there exists unique $j\in\{0,1,\ld,i_0\}$ such that
\begin{eqnarray}
|\la^{\mi_{v_2}+j}\bi_{j+1}-1|<\e_1.  \label{eq3}
\end{eqnarray}
We remark by the previous Lemma \ref{lem3.4} that for every $\la\in[a_0,b_0]$ there exists unique $a_i$, $i\in\{0,1,\ld,i_0\}$ such that $a_i\le\la\le a_{i+1}$ if $i\le i_0-1$ or $a_{i_0}\le\la\le b_0$, where the numbers $a_i$, $i=0,1,\ld,i_0$ are defined completely by Lemma \ref{lem3.4}.

So, we have defined completely the positive numbers $\bi_{j+1}$, for $j=0,1,\ld,i_0$.

Now, we define $x_{\mi_{v_2+i}+j}=\bi_{i+1}q_j$ for every $i=0,1,\ld,i_0$ and for every $j=1,2,\ld,v_0$.

The previous terms $x_{\mi_{v_2+i}+j}$ for $i=0,1,\ld,i_0$, $j=1,2,\ld,v_0$ are defined well because $\mi_{n+1}-\mi_n>v_0$ for every $n=1,2,\ld$ by the definition of the sequence $(\mi_n)$.

Finally, we define that $x_j=0$ for every $j\in\N$, $j>\mi_{v_2}$ for which there exists not $i\in\{0,1,\ld,i_0\}$ and $j_1\in\{1,2,\ld,v_0\}$ such that $j=\mi_{v_2+i}+j_1$.

By the previous procedure we have defined completely the vector\\ $x_0=(x_1,x_2,\ld,x_n,\ld)$ where $\{n\in\N\big|x_n\neq0\}$ is finite and thus $x_0\in\el^2$ obviously.

Now, we show that the vector $x_0$ satisfies relation (\ref{eq1}).

Firstly we prove that $x_0\in S_{\el^2}(c_0,\e_0)$.

By the definition of the vector $x_0$ we get:
\begin{align}
\|x_0-c_0\|^2:&=\sum^{+\infty}_{j=1}|x_j-c_j|^2=\sum^{+\infty}_{j=v_1+1}|x_j|^2 \nonumber\\
&=\sum^{v_0}_{j=1}\sum^{i_0}_{i=0}|x_{\mi_{v_2+i}+j}|^2\nonumber\\
&=\sum^{v_0}_{j=1}\sum^{i_0}_{i=0}|\bi_{i+1}q_j|^2=\sum^{v_0}_{j=1}|q_j|^2
\sum^{i_0}_{i=0}|\bi_{i+1}|^2 \nonumber \\
&\le v_0M^2_1\cdot\sum^{i_0}_{i=0}|\bi_{i+1}|^2=v_0M^2_1\cdot\sum^{i_0}_{i=0}
\frac{1}{a_i^{2\mi_{v_2+i}}}\nonumber\\
&\le v_0M^2_1\cdot\sum^{i_0}_{i=0}\frac{1}{a_0^{2\mi_{v_2+i}}}<v_0M^2_1\cdot
\sum^{+\infty}_{v=2\mi_{v_2}}\frac{1}{a^v_0} \nonumber\\
&=v_0M^2_1\frac{1}{a_0^{2\mi_{v_2}}}\cdot\frac{1}{1-\dfrac{1}{a_0}}<\e^2_0 \label{eq4}
\end{align}
by the inequality (c) above for $\mi_{v_2}$.

Inequality (\ref{eq4}) gives that $x_0\in S_{\el^2}(c_0,\e_0)$, so property (i) of (\ref{eq1}) holds.

We show now that property (ii) also holds.

We fix some $\la\in[a_0,b_0]$. Then there exists unique $\rho_0\in\{1,2,\ld,i_0-1\}$ such that $a_{\rho_0}\le\la<a_{\rho_0}\cdot(1+\e_0)^{\frac{1}{\mi_{v_2}+\rho_0}}$ or $a_{i_0}\le\la\le b_0$.

We show that
\[
\|(\la B)^{\mi_{v_2+\rho_0}}(x_0)-y_{j_0}\|<\frac{1}{s_0}.
\]
We have:
\begin{align}
\|(\la B)^{\mi_{v_2+\rho_0}}(x_0)-y_{j_0}\|^2=&\sum^{v_0}_{j=1}|\la^{\mi_{v_0+\rho_0}}
x_{\mi_{v_2+\rho_0}+j}-q_j|^2 \nonumber \\
&+\sum^{+\infty}_{j=v_0+1}|\la^{\mi_{v_2+\rho_0}}
x_{\mi_{v_2+\rho_0}+j}|^2.  \label{eq5}
\end{align}
By definition we have for $j=1,2,\ld,v_0$\; $x_{\mi_{v_2+\rho_0}+j}=\bi_{\rho_0+1}q_j$. So, for $j=1,2,\ld,v_0$ and
$\la\in[a_{\rho_0},a_{\rho_0+1}]$, where
$ a_{\rho_0+1}=a_{\rho_0}\cdot(1+\e_1)^{\frac{1}{\mi_{v_2}+\rho_0}}$
or $\la\in[a_{\rho_0},b_0]$ if $\rho_0=i_0$ we have
\begin{align*}
|\la^{\mi_{v_2+\rho_0}}x_{\mi_{v_2+\rho_0}+j}-q_j|^2&=|\la^{\mi_{v_2+\rho_0}}
\bi_{\rho_0+1}q_j-q_j|^2 \\
&=|\la^{\mi_{v_2+\rho_0}}\bi_{\rho_0+1}-1|^2|q_j|^2 \\
&\le|\la^{\mi_{v_2+\rho_0}}\bi_{\rho_0+1}-1|^2\cdot M^2_1<\e^2_1M^2_1=\frac{1}{2v_0s^2_0}.
\end{align*}
So we have:
\begin{eqnarray}
\sum^{v_0}_{j=1}|\la^{\mi_{v_2+\rho_0}}x_{\mi_{v_2+\rho_0}+j}-q_j|^2<\frac{1}{2s^2_0}. \label{eq6}
\end{eqnarray}
If $\rho_0=i_0$ the second member of (\ref{eq5}) is 0 and the conclusion holds by (\ref{eq6}).

So for the sequel we suppose that $\rho_0\le i_0-1$. In this case we get
\begin{align*}
\sum^{+\infty}_{j=v_0+1}|\la^{\mi_{v_2}+\rho_0}x_{\mi_{v_2+\rho_0}+j}|^2&=\la^{2\mi_{v_2+\rho_0}}
\sum^{v_0}_{j=1}\sum^{i_0-\rho_0}_{i=1}|x_{\mi_{v_2+\rho_0+i}+j}|^2 \nonumber\\
&=\la^{2\mi_{v_2+\rho_0}}\sum^{v_0}_{j=1}|q_j|^2\sum^{i_0-\rho_0}_{i=1}|\bi_{\rho_0+i+1}|^2 \nonumber \\
&\le\la^{2\mi_{v_2+\rho_0}}v_0M^2_1\sum^{i_0-\rho_0}_{i=1}|\bi_{\rho_0+i+1}|^2 \nonumber \\
&=\la^{2\mi_{v_2+\rho_0}}v_0M^2_1\sum^{i_0-\rho_0}_{i=1}\frac{1}{\Big(a^{\mi_{v_2+\rho_0+i}}_{\rho_0+i}\Big)^2} \nonumber \\
&\le\la^{2\mi_{v_2+\rho_0}}v_0M^2_1\sum^{i_0-\rho_0}_{i=1}\frac{1}{a^{2\mi_{v_2+\rho_0+i}}_{\rho_0+1}} \nonumber \\ &<\la^{2\mi_{v_2+\rho_0}}v_0M^2_1\cdot\sum^{+\infty}_{v=2\mi_{v_2+\rho_0+1}}
\frac{1}{a^v_{\rho_0+1}}\nonumber \\
&=\la^{2\mi_{v_2+\rho_0}}v_0M^2_1\frac{1}{a_{\rho_0+1}^{2\mi_{v_2+\rho_0+1}}}\cdot
\frac{1}{1-\dfrac{1}{a_{\rho_0+1}}}\nonumber
\end{align*}
\begin{align}
&<a^{2\mi_{v_2+\rho}}_{\rho+1}v_0M^2_1
\frac{1}{a^{2\mi_{v_2+\rho_0+1}}_{\rho_0+1}}\cdot\frac{1}{1-\dfrac{1}{a_0}}\nonumber \\
&=\frac{v_0M^2_1}{1-\dfrac{1}{a_0}}\cdot\frac{1}{a_{\rho_0+1}^{2(\mi_{v_2+\rho_0+1}-
\mi_{v_2+\rho_0})}} \nonumber \\
&<\frac{v_0M^2_1}{1-\dfrac{1}{a_0}}\cdot\frac{1}{a_0^{2(\mi_{v_2+\rho_0+1}-\mi_{v_2+\rho_0})}}<
\frac{1}{2s^2_0}  \label{eq7}
\end{align}
because $\mi_{v+1}-\mi_v>M_2$ for every $v\ge v_2$ by the hypothesis (i) for the sequence $(\mi_n)$.

By (\ref{eq5}), (\ref{eq6}) and (\ref{eq7}) we get that $\|(\la B)^{\mi_{v_2+\rho_0}}(x_0)-y_{j_0}\|<\dfrac{1}{s_0}$ for the arbitrary $\la\in[a_0,b_0]$. This completes property (ii) of (\ref{eq1}) and the proof of this Lemma \ref{lem3.6} is completed. \qb
\end{Proof}

Now by Lemmas \ref{lem3.2}, \ref{lem3.3}, \ref{lem3.6} and the facts that the space $(\el^2,\|\cdot\|)$ is a complete metric space and Baire's category Theorem the proof of Theorem \ref{thm3.1} is completed.
\section{A result in measure and category}\label{sec4}
\noindent
In this section we prove item $(iii)$ of Theorem \ref{mainresult}. Actually, we shall prove the following, more general, result. As we already mentioned in the Introduction, its proof elaborates on the proof of Proposition 5.2 from \cite{BaCo}.

\begin{thm}\label{thm4.1}
Let $(k_n)$ be a strictly increasing sequence of positive integers. Let $T$ be a bounded linear operator acting on a (complex) Banach space $X$ such that $((\la T)^{k_n} )$ is hypercyclic for every $\la \in \mathbb{C}$ with $|\la |>1$ and assume in addition that 
$$\bigcap\limits_{\la\in (1,+\infty)}\ch C(\{(\la T)^{k_n}\})\neq \emptyset .$$ Then, there exists a $G_{\delta }$ and dense subset $\cp$ in $\{ \la \in \mathbb{C}: |\la |>1 \}$ with full $2$-dimensional Lebesgue measure in $ \{ \lambda \in \mathbb{C}: |\la |>1 \}$ such that $\bigcap_{\la\in\cp}\ch C(\{(\la T)^{k_n}\}) $ is  residual in $X$. In particular, $\bigcap_{\la\in\cp}\ch C(\{(\la T)^{k_n}\})\neq\emptyset$. 
\end{thm}
\begin{Proof}
By performing a change of variables it suffices to prove the following: \\
{\bf Claim}. \textit{Fix $x\in \bigcap\limits_{\la\in (1,+\infty)}\ch C(\{(\la T)^{k_n}\})$. Then there exists a $G_{\de }$ and dense subset $A$ of $(1,+\infty )\times \mathbb{R}$ with full ($2$-dimensional) Lebesgue measure such that the set $\{ ((rT)^{k_n}x, e^{2\pi ik_n\theta }) :n=1,2,\ldots \}$ is dense in $X\times \mathbb{T}$} for every $(r,\theta ) \in A$.

Here $\mathbb{T}$ denotes the unit circle, i.e. $\mathbb{T}=\{ z\in \mathbb{C}:|z|=1 \}$.\\
{\bf Proof of Claim}. Let $\{ x_j: j\in \mathbb{N} \}$, $\{ t_l:l\in \mathbb{N} \}$  be dense subsets of $X$, $\mathbb{T}$ respectively. For every $j,l,s,n \in \mathbb{N}$ define the set
$$A_{j,l,s,n}:=\left\{ (r,\theta )\in (1,+\infty )\times \mathbb{R}: \| (rT)^{k_n}x-x_j\| <\frac{1}{s}, |e^{2\pi ik_n\theta }-t_l|<\frac{1}{s} \right\}.$$ We shall prove that the set $A:=\bigcap_{j,l,s}\bigcup_nA_{j,l,s,n}$ has the desired properties. Since $A_{j,l,s,n}$ is open we conclude that $A$ is $G_{\de}$. Let us show that $A$ is dense in $(1,+\infty )\times \mathbb{R}$. In view of Baire's theorem it suffices to prove that for any fixed $j,l,s\in \mathbb{N}$ the set $\bigcup_{n}A_{j,l,s,n}$ is dense in $(1,+\infty )\times \mathbb{R}$. To this end, fix $j,l,s\in \mathbb{N}$ and let $b>1$, $a\in \mathbb{R}$ and $\epsilon >0$. We seek $r>0$, $\theta \in \mathbb{R}$ and $n\in \mathbb{N}$ such that 
$$|b-r|<\epsilon, \,\,|a-\theta |<\epsilon ,\,\, |t_l-e^{2\pi i k_n\theta}| \,\, \textrm{and}\,\, \| (rT)^{k_n}x-x_j\| <1/s.$$ 
Define the set $B:=\{ k_n: \| (bT)^{k_n}x-x_j\| <1/s \}$ and consider its elements in an increasing order, say $k_{\rho_1}<k_{\rho_2}<\cdots $. Of course, we have $B=\{ k_{\rho_n}: n\in \mathbb{N} \}$. Now we use Weyl's theorem, see Theorem 4.1 in page $32$ from \cite{KN}, to conclude that the sequence $(k_{\rho_n}\theta )$ is uniformly distributed modulo $1$ for almost all $\theta $ in $\mathbb{R}$. Hence, there exists $\theta \in \mathbb{R}$ such that the set $\{ e^{2\pi ik_{\rho_n}\theta } :n\in \mathbb{N} \}$ is dense in $\mathbb{T}$ and $|a-\theta |<\epsilon $. Finally, setting $r:=b$ and from all the above we conclude that there exists $n:=\rho_m$ for some $m\in \mathbb{N}$ such that    
$$|b-r|<\epsilon, \,\,|a-\theta |<\epsilon ,\,\, |t_l-e^{2\pi i k_n\theta}| \,\, \textrm{and}\,\, \| (rT)^{k_n}x-x_j\| <1/s, $$ which is what we wanted to prove. 
It remains to show that $A$ has full measure in $(1,+\infty )\times \mathbb{R}$. Actually, it is enough to prove that the set $\bigcup_{n}A_{j,l,s,n}$ has full measure in $(1,+\infty )\times \mathbb{R}$ for every $j,l,s\in \mathbb{N}$. Fix $j,l,s\in \mathbb{N}$ and take any four numbers $d_1,d_2,d_3,d_4$ with $d_1<d_2$, $1<d_3<d_4$. For any subset $B$ of $(1,+\infty )\times \mathbb{R}$ the symbol $B_r$ stands for its section, i.e. $B_r:=\{ \theta \in \mathbb{R}: (r,\theta )\in E \}$ and for simplicity reasons we set $E:=\bigcup_nA_{j,l,s,n}$. Observe that the proof of denseness result implies that for every $r\in [d_3,d_4]$ 
we have $(r,\theta )\in E$ for almost every $\theta $ in $\mathbb{R}$ (of course the set of such $\theta$'s depends on $r$). It now follows that
$$ \mu ((E\cap ([d_3,d_4]\times [d_1,d_2]))_r)=d_2-d_1=\mu (( [d_3,d_4]\times [d_1,d_2])_r) \,\, \textrm{for every} \,\, r\in [d_3,d_4] ,$$ where $\mu$ denotes the Lebesgue measure,
and by Cavalieri's principle, see page 149 in \cite{Hal}, we conclude that 
$$\mu \times \mu (E\cap ([d_3,d_4]\times [d_1,d_2]))=(d_2-d_1)(d_4-d_3).$$ 
Thus, $E$ has full measure in $(1,+\infty )\times \mathbb{R}$. This completes the proof of the Claim and hence that of Theorem \ref{thm4.1}.
\end{Proof}
Item $(ii)$ of Theorem \ref{mainresult} and Theorem \ref{thm4.1} directly imply item $(iii)$ of Theorem \ref{mainresult}.

Department of Mathematics and Applied Mathematics, University of Crete, Panepistimiopolis Voutes, 700-13, Heraklion, Crete, Greece.\\
email:tsirivas@uoc.gr


\begin{thebibliography}{99}
\bibitem{AbGo} E. Abakumov, J. Gordon, Common hypercyclic vectors for multiples of backward shift,
J. Funct. Anal. 200 (2003), 494-504.
%
\bibitem{Ba} F. Bayart, Common hypercyclic vectors for high dimensional families of operators, arXiv:1503.08574 

\bibitem{BaCo} F. Bayart and G. Costakis, Hypercyclic operators and rotated orbits with polynomial phases, Journal of the London Mathematical Society 89 (2014), 663-679.

\bibitem{2}  F. Bayart and E. Matheron, Dynamics of linear operators, Cambridge Tracts in Math. 179, Cambridge Univ. Press, 2009.
%
\bibitem{3}  F. Bayart and E. Matheron, How to get common universal vectors, Indiana Univ. Math. J., 56 (2007), 553-580.
%

\bibitem{ChaSa2} K. C. Chan, R. Sanders, Two criteria for a path of operators to have common hypercyclic vectors, J.
Operator Theory 61 (2009), 191-223.

\bibitem{ChaSa4} K. C. Chan, R. Sanders, An SOT-dense path of chaotic operators with same hypercyclic vectors, J. Operator Theory 66 (2011), 107-124.

\bibitem{5} A. Conejero, V. M\"{u}ller, A. Peris, Hypercyclic behaviour of operators in a hypercyclic $C_0$-semigroup, J. Funct. Anal., 244 (2007), 342-348.
%
%
\bibitem{CoSa} G. Costakis and M. Sambarino, Genericity of wild holomorphic functions and common hypercyclic vectors, Adv. Math. 182 (2004), 278-306.

\bibitem{8} G. Costakis, N. Tsirivas and V. Vlachou, Non-existence of common hypercyclic vectors for certain families of translations operators, to appear in Computational Methods and Function Theory. 
%
\bibitem{10} K. Grosse-Erdmann and A. Peris, Linear Chaos, Universitext, 2011, Springer.
%
\bibitem{Hal} P. R. Halmos, Measure Theory, Graduate Texts in Math., 1974 Springer-Verlag. 

\bibitem{KN} L. Kuipers and H. Niederreiter, Uniform distribution of sequences, 2006, Dover Publications.

\bibitem{LeMu} F. Leon-Saavedra, Fernando, V. M\"{u}ller, Rotations of hypercyclic and supercyclic operators,
Integral Equations Operator Theory 50 (2004), 385-391.
%


%
\bibitem{San} R. Sanders, Common hypercyclic vectors and the hypercyclicity criterion, Integral Equations Operator
Theory 65 (2009), 131-149.

\bibitem{Shka} S. Shkarin, Universal elements for non-linear operators and their applications,
J. Math. Anal. Appl. 348 (2008), 193-210.

\bibitem{13} S. Shkarin, Remarks on common hypercyclic vectors, J. Funct. Anal. 258, (2010), 132-160.

\bibitem{Tsi} N. Tsirivas, Existence of common hypercyclic vectors for translation operators,  arXiv:1411.7815

\bibitem{Tsi2} N. Tsirivas, Common hypercyclic functions for translation operators with large gaps,  arXiv:1412.0827

\bibitem{Tsi3} N. Tsirivas, Common hypercyclic functions for translation operators with large gaps II,  arXiv:1412.1963






\end{thebibliography}
\end{document}